\def\ZZ         {{\mathbb Z}}
\def\RR         {{\mathbb R}}
\def\CC         {{\mathbb C}}
  
\def\HH         {{\mathbb H}}
\def\QQ         {{\mathbb Q}}
\def\PP         {{\mathbb P}}

\def\E           {{\cal E}}
\def\K           {{\cal K}}

\def\ii         {{\rm i}}
\def\ee         {{\rm e}}

\def\dim        {{\rm dim}}

\def\sin        {{\rm sin}}

\def\sinh       {{\rm sinh}}

\def\Ell        {{\cal ELL}}

\def\Prod       {{\displaystyle\prod}}

\def\cal        {\mathcal}

\documentclass{amsart}
\usepackage{amssymb}
\usepackage{verbatim}
\usepackage{eucal}
\usepackage{mathrsfs}
\usepackage{graphicx}
\usepackage{psfrag}

\newtheorem{theorem}{Theorem}[section]
\newtheorem{lemma}[theorem]{Lemma}
\newtheorem{prop}[theorem]{Proposition}

\theoremstyle{definition}
\newtheorem{dfn}[theorem]{Definition}

\newtheorem{example}[theorem]{Example}

\theoremstyle{remark}
\newtheorem{remark}[theorem]{Remark}

\begin{document}
  
\title{Elliptic genera,
real algebraic varieties and quasi-Jacobi forms}

\author{Anatoly Libgober}\footnote{author was partially supported by NSF grant}
\address{Department of Mathematics\\
University of Illinois\\
Chicago, IL 60607}
\email{libgober@math.uic.edu}

\begin{abstract} This paper surveys the push forward formula for 
elliptic class and various applications obtained in the papers by 
L.Borisov and the author. In the remaining part we discuss 
the ring of quasi-Jacobi forms which allows to characterize 
the functions which are the elliptic genera of almost complex manifolds
and extension of Ochanine elliptic genus to 
certain singular real algebraic varieties.
\end{abstract}

\maketitle

\section*{Introduction}

Interest in elliptic genus of complex manifolds stems from its 
appearance in a wide 
variety of geometric and topological problems.
At the first glance, this is an invariant of
complex cobordism class modulo torsion and hence depends only on the 
Chern numbers
of manifold. On the other hand elliptic genus is a 
holomorphic function of $\CC \times \HH$ where 
$\HH$ is the upper half plane. In one of heuristic
approaches, elliptic genus 
is an index of an operator on the loop space (\cite{witten87}) 
and as such 
it has counterparts defined for $C^{\infty}$, oriented or Spin manifolds
(studied prior to the study of complex case 
(cf. \cite{Ochanine}).
It comes up in the study of geometry and topology of loop spaces and, 
more specifically,  the 
chiral deRham complex (cf. \cite{MSV}), in the study of 
invariants of singular algebraic 
varieties (\cite{Duke}), in particular orbifolds,
 and more recently in the study of 
Gopakumar-Vafa 
and Nekrasov conjectures (cf. \cite{liu} \cite{GL}). 
It is closely related to the fast 
developing subject of elliptic cohomology (cf. \cite{Thomas}). 
There are various versions
of elliptic genus: equivariant, higher elliptic genus
obtained by twisting by cohomology classes of the fundamental
group etc, elliptic genus of pairs and orbifold elliptic genus. 
There is interesting connection with singularities
of weighted homogeneous polynomials (so called Landau-Ginzburg models).

In present note we shall review several recent developments in 
the study of elliptic genus and refer a reader to review 
\cite{oldreview} for additional details on earlier results.
Then we shall focus on its aspects of elliptic genus: 
its extension to real singular varieties and its modularity property
(or the lack of it). Extension of elliptic genus to real 
singular varieties was suggested by B.Totaro (\cite{Totaro}) and 
our approach is based 
on the push forward formula for the elliptic class 
used in \cite{Duke} to extend elliptic genus from smooth to certain 
singular complex projective 
varieties. 

In section \ref{start} we shall discuss this push 
forward formula. It appears as the main technical tool
in many applications mentioned later in the paper.
The rest of section \ref{ellipticgenus}  
discusses the relation with other invariants and series of 
applications based on the material of works \cite{annals},\cite{Duke}
and \cite{invent}. It includes a discussion of a relation between elliptic
genus and $E$-function, application to McKay correspondence, elliptic 
genera of non-simply-connected manifolds (higher elliptic genera) 
and generalizations of a formula of R.Dijkgraaf, Moore, E.Verlinde, 
H.Verlinde. 
Other applications 
in equivariant context are discussed in R.Waelder paper in this volume.
The proof of independence of resolutions of our definition of elliptic genus 
for certain real algebraic varieties is
 given in section \ref{lastsection}. 

Section 2 deals with modularity properties of elliptic genus.
In the Calabi Yau cases (of pairs, orbifolds etc.) elliptic genus is 
a weak Jacobi form (cf. definition below).
Also it is important to have description of the elliptic genus 
in non Calabi Yau situations not just as 
a function on $\CC \times \HH$ 
but as an element of a finite dimensional algebra of functions.
It turns out that in 
the absence of Calabi Yau condition the elliptic genus belong to 
a very interesting algebra of functions on $\CC \times \HH$, 
which we call the algebra of quasi-Jacobi forms, 
and which is only slightly bigger than the algebra of weak Jacobi forms. 
This algebra of quasi-Jacobi forms 
is a counterpart of quasi-modular forms (cf. \cite{Kaneko})
and is related to elliptic 
genus in the same way as quasi-modular forms are related to the Witten 
genus (cf. \cite{ZagierLandweb}). The algebra of quasi-Jacobi forms 
is generated by certain two variable Eisenstein series (masterfully reviewed 
by A.Weil in  \cite{Weil}) and has many properties 
parallel to the properties in quasi-modular case. A detailed description 
of properties of quasi-Jacobi forms appears to be absent in the literature and 
we discuss the algebra of these forms in section 2 (for more
details see its introduction). We end this section with discussion 
of differential operators Rankin-Cohen brackets on the space of Jacobi forms.

Finally in section \ref{sectionrealsing} construct extension of Ochanine
genus to real algebraic varieties with certain class of singularities.
This extends some results of B.Totaro in \cite{Totaro}.

For readers convenience we give ample references to 
prior work on elliptic genus where more detailed 
information on the subject can be obtained.
Part 2 dealing with quasi-Jacobi forms can be read 
independently of the rest of the paper.

The author wants to express his gratitude to Lev Borisov. 
The material in section 1 is a survey of joint papers with him
and results of section 2 are based on discussions 
with him several years ago.

\section{Elliptic genus.}\label{ellipticgenus}

\subsection{Elliptic genus of singular varieties and push-forward
 formulas}\label{start}
Let $X$ be a projective manifold. We shall use the Chow groups $A_*(X)$
with complex coefficients (cf. \cite{Fulton}). Let $F$ the ring of functions 
on $\CC \times \HH$ where $\HH$ is the upper half-plane.
The elliptic class of $X$ is an element in $A_*(X) \otimes_{\CC} F$ 
given by:
\begin{equation}\label{ellipticclass}
\Ell(X)=\prod_i x_i {{\theta ({{x_i} \over {2 \pi \ii}}-z,\tau)} \over
{\theta ({{x_i} \over {2 \pi \ii }}, \tau)}}[X] 
\end{equation}
where 
\begin{equation}
\theta(z,\tau)=q^{1 \over 8}  (2 \sin \pi z)
\prod_{l=1}^{l=\infty}(1-q^l)
 \prod_{l=1}^{l=\infty}(1-q^l \ee^{2 \pi \ii z})(1-q^l \ee^{-2 \pi \ii
z})
\end{equation}
is the Jacobi theta function considered as an element 
in $F$ with $q=e^{2 \pi i \tau}$ 
(cf. \cite{Chandra}), $x_i$ are the Chern roots of the tangent bundle of $X$ 
and $[X]$ is the fundamental class of $X$. The component $Ell(X)$ 
in $A_0(X)=F$ is the elliptic genus of $X$. 

Components of each degree of (\ref{ellipticclass}),
evaluated on a class in $A^*(X)$,
are linear combination of symmetric functions in $c_i$ i.e. the Chern classes
of $X$. In particular $Ell(X)$ depends only on the class of $X$ in the 
ring $\Omega^U \otimes \QQ$ of unitary cobordisms. 

The homomorphism 
$\Omega^U\otimes \QQ \rightarrow F$ taking $X$ to $Ell(X)$ can be described
without reference to theta functions as in (\ref{ellipticclass}).
Let $M_{A_1,3}$ be the class of complex analytic space 
``having only $A_1$-singularities in codimension 3'' i.e. 
having the singularities of the following type: the singular set $Sing X$ of 
$X \in M_{A_1,3}$ is {\it a manifold} such that 
${\rm dim}_{\CC} Sing X={\rm dim} X-3$ and for an embedding $X \rightarrow 
Y$ where $Y$ is a manifold  and a transversal $H$ to $Sing X$ in $Y$, 
the pair $(H \cap X, H \cap Sing X)$ is analytically equivalent to 
the pair $(\CC^4,H_0)$ where $H_0$ is given by $x^2+y^2+z^2+w^2=0$.
Each $X \in M_{A_1,3}$ admits two small resolutions $\tilde X_1 \rightarrow X$
and  $\tilde X_2 \rightarrow X$
in which the exceptional set is fibration over $Sing X$ with the fiber 
$\PP^1$. One says that manifolds underlying the resolutions obtained
from each other by classical flop.

\begin{theorem}\label{totarotheorem}(cf. \cite{Totaroannals}) The kernel of 
the homomorphism $Ell: \Omega^U\otimes \QQ
 \rightarrow F$ taking an almost complex
manifold $X$ to its elliptic genus $Ell(X)$ is 
the ideal generated by classes of differences
$\tilde X_1-\tilde X_2$ of two small resolutions of a variety 
in $M_{A_1,3}$. 
\end{theorem}

More generally one can fix a class of singular spaces and a type of 
resolutions and consider the quotient of $\Omega^U \otimes \QQ$ by the ideal 
generated by differences of manifolds underlying resolutions of the 
same analytic space. The quotient map by this ideal $\Omega^U 
\otimes \QQ \rightarrow R$
provides a genus and hence the collection of Chern numbers (linear
combination of Chern monomials $c_{i_1} \cdot ....\cdot c_{i_k}[X]$ 
($\sum i_s={\rm dim} X$) which can be made explicit via Hirzebruch 
procedure with generating series (cf. \cite{hirztop}). 
These are the Chern numbers which can 
be defined for the chosen class of singular varieties and chosen 
class of resolutions. The ideal in theorem \ref{totarotheorem}, 
it turns out, 
corresponds to a much larger classes of singular spaces and resolutions.
This method of defining Chern classes of singular varieties is 
an extension of the philosophy underlying the question of 
Goresky and McPherson (cf. \cite{Goreski}): 
which Chern numbers can be defined via 
resolutions independently of resolution.

\begin{dfn} An analytic space 
$X$ is called $\bf Q$-Gorenstein if the divisor 
$D$ of a meromorphic form $df_1 \wedge .... \wedge df_{{\rm dim}X}$
is such that for some $n \in \ZZ$ the divisor $nD$ in locally principal
(i.e. $K_X$ is $\QQ$-Cartier). In particular for any codimension one component 
$E$ of the exceptional divisor of a map $\pi: \tilde X \rightarrow X$ 
the multiplicity $a_E={\rm mult}_E \pi^*(K_X)$ is well defined 
and a singularity is called log-terminal if there is 
resolution $\pi$ such that $K_{\tilde X}=\pi^*(K_X)+\sum a_E E$
and $a_E >-1$. A resolution is called crepant if $a_E=0$.
\end{dfn}

\begin{theorem}\label{duketheorem}(\cite{Duke}) 
The kernel of the elliptic genus $\Omega^U \otimes \QQ \rightarrow F$ 
is generated by the differences of $\tilde X_1-\tilde X_2$ 
of manifolds underlying crepant resolutions of the singular spaces
with $\QQ$-Gorenstein singularities admitting crepant resolutions.
\end{theorem}

The proof of theorem (\ref{duketheorem}) is based on extension of 
elliptic genus $Ell(X)$ of manifolds to the elliptic genus of 
pairs $Ell(X,D)$ where $D$ is a divisor on $X$ having the normal 
crossings as the only singularities. This is similar to the 
situation in the study of motivic $E$-functions of quasi-projective 
varieties (cf.\cite{bat},\cite{looijenga}).   
In fact, motivation from 
other problems (e.g. study of McKay correspondence cf. \cite{annals}) 
suggests looking at the triples $(X,D,G)$ where $X$ is a normal 
variety, $G$ is a finite group acting on $X$ and to introduce the elliptic
class $\Ell(X,D,G)$ (cf. \cite{annals}). More precisely, let 
$D$ be a $\QQ$-divisor i.e. $D=\sum a_iD_i$ with 
$D_i$ being irreducible and $a_i \in \QQ$.
The pair $(X,D)$ is called Kawamata log-terminal 
(cf. \cite{KMM}) if $K_X+D$ is $\QQ$-Cartier and there
is a birational morphism $f: Y \rightarrow X$ where $Y$ is 
smooth and the union of 
the proper preimages of components of $D$ and the components of 
exceptional set  $E=\bigcup_{\cal J} E_j$
 form a normal crossing divisor  
such that $K_Y=f^*(K_X+\sum a_iD_i)+\sum \alpha_jE_j$ 
where $\alpha_j >-1$
(here 
$K_X,K_Y$ are the canonical classes of $X$ and $Y$ respectively).
The triple $(X,D,G)$ where $X$ is a non-singular variety, 
$D$ is a divisor and $G$ is a finite group of biholomorphic automorphisms is  
 called $G$-normal (cf. \cite{bat}, \cite{annals}) 
if the components of $D$ form a normal crossings divisor 
and the isotropy group of any point acts trivially on the 
components of $D$ containing this point.

\begin{dfn}(cf. \cite{annals} definition 3.2)\label{orbifoldelliptic}
 Let $(X,E)$ be a Kawamata log terminal $G$-normal pair (i.e. in 
particular $X$ is smooth with $D$ being a normal crossing divisor)
 and $E=-\sum_{k \in {\cal K}} \delta_kE_k$. 
The {\it orbifold} elliptic class of 
 $(X,E,G)$ is the class in $A_*(X,\QQ)$ given by:

\begin{equation}\label{formulaorbifoldelliptic}
{\cal Ell}_{orb}(X,E,G;z,\tau):=
\frac 1{\vert G\vert }\sum_{g,h,gh=hg}\sum_{X^{g,h}}[X^{g,h}]
\Bigl(
\prod_{\lambda(g)=\lambda(h)=0} x_{\lambda} 
\Bigr)
\end{equation}
$$\times\prod_{\lambda} \frac{ \theta(\frac{x_{\lambda}}{2 \pi \ii }+
 \lambda (g)-\tau \lambda(h)-z )} 
{ \theta(\frac{x_{\lambda}}{2 \pi \ii }+
 \lambda (g)-\tau \lambda(h))}  \ee^{2 \pi \ii \lambda(h)z}
$$
$$\times\prod_{k}
\frac
{\theta(\frac {e_k}{2\pi\ii}+\epsilon_k(g)-\epsilon_k(h)\tau-(\delta_k+1)z)}
{\theta(\frac {e_k}{2\pi\ii}+\epsilon_k(g)-\epsilon_k(h)\tau-z)}
{}
\frac{\theta(-z)}{\theta(-(\delta_k+1)z)} \ee^{2\pi\ii\delta_k\epsilon_k(h)z}.
$$
\end{dfn}
\noindent 
where $X^{g,h}$ denotes an irreducible component of the fixed set of the 
commuting elements $g$ and $h$ and $[X^{g,h}]$ denotes the image 
of the fundamental class in $A_*(X)$. The restriction of $TX$ to $X^{g,h}$
splits into linearized bundles according to the ($[0,1)$-valued) characters 
$\lambda$  of $\langle g,h\rangle$ (sometimes denoted $\lambda_W$ where 
$W$ is a component of the fixed point set). 
Moreover, $e_k=c_1(E_k)$ 
and $\epsilon_k$ is the character of ${\cal O}(E_k)$ restricted to $X^{g,h}$ 
if $E_k$ contains $X^{g,h}$ and is zero otherwise.

One would like to define the elliptic genus of a 
Kawamata log-terminal pair $(X_0,D_0)$ as 
(\ref{formulaorbifoldelliptic}) calculated 
for a $G$-equivariant resolution $(X,E) \rightarrow (X_0,D_0)$.
Independence of (\ref{formulaorbifoldelliptic}) of
resolution  and the proof of 
(\ref{duketheorem}) both depend on the following push-forward formula:

\begin{theorem}\label{orbblowup}
Let $(X,E)$ be a Kawamata log-terminal $G$-normal pair
and let $Z$ be a smooth $G$-equivariant locus in
$X$ which is normal crossing to $Supp(E)$.
Let $f:\hat X\to X$ denote the blowup of $X$ along $Z$. We define 
$\hat E$ by $\hat E=-\sum_{k}\delta_k\hat E_k-\delta Exc(f)$
where $\hat E_k$ is the proper transform of $E_k$ and 
$\delta$ is determined from $K_{\hat X}+\hat E = f^*(K_X + E)$.
Then $(\hat X,\hat E)$ is a Kawamata log-terminal $G$-normal pair and 
\begin{equation}\label{orbblowupformula}
f_*{\cal Ell}_{orb}(\hat X,\hat E,G;z,\tau)={\cal Ell}_{orb}(X,E, G;z,\tau).
\end{equation}
\end{theorem}

Independence of resolution is a consequence of the weak factorization 
theorem (cf. \cite{weak})
and (\ref{orbblowup}), while theorem (\ref{duketheorem}) 
follows since both 
${\cal Ell}(X_1)$ and ${\cal Ell}(X_2)$ coincide with the elliptic 
genus of the pair $(\tilde X,\tilde D)$ where $\tilde X$ is a resolution 
of $X$ dominating both $X_1$ and $X_2$ (where $D=K_{\tilde X/X}$ cf.
\cite{Duke} Prop. 3.5 and also \cite{Wang}). For discussion of 
orbifold elliptic genus on orbifolds more general than just global quotients 
see \cite{Dong}.

\subsection{Relation to other invariants.}\label{relationtoother}
 V.Batyrev in \cite{bat}
for a $G$-normal triple $(X,D,G)$ defined $E$-function
$E_{orb}(X,D,G)$ depending on the Hodge theoretical invariants
(there is also motivic version cf.\cite{bat},\cite{looijenga}).
Firstly for a quasi-projective algebraic variety $W$ one put
(cf. \cite{bat}, Definition 2.10):
\begin{equation}\label{efunction}
E(W,u,v)=(-1)^i\sum_{p,q} \dim Gr^p_F Gr^{p+q}_W(H^i_c(W,\CC))u^pv^q
\end{equation}
where $F$ and $W$ are Hodge and weight filtrations of Deligne's mixed
Hodge structure (\cite{deligne}).
In particular $E(W,1,1)$ is the topological Euler characteristic 
of $W$ (with compact support). If $W$ is compact then one obtains 
Hirzebruch's $\chi_y$-genus:
\begin{equation}\label{chisuby} 
\chi_y(W)=\sum_{i,j} (-1)^q\dim H^q(\Omega^p_W)y^p
\end{equation} 
(cf. \cite{hirztop}) for $v=-1,u=y$ and hence the arithmetic genus, 
signature etc. are special values of (\ref{efunction}).
Secondly, for a $G$-normal pair as in (\ref{orbifoldelliptic}) one 
stratifies $D=\bigcup_{k \in K} D_k$ 
by strata 
$D^{\circ}_{J}=
\bigcap_{j \in J} D_j-\bigcup_{k \in K-J} D_k, \ (J \subset K,
 \ \bigcap_{j \in J} D_j=X \ {\rm for}\  J=\emptyset)$
 and defines
\begin{equation}\label{efunctionorb}
E(X,D,G,u,v)=
\end{equation}
$$\sum_{\{g\},W \subset X^g} 
(uv)^{\sum \epsilon_{D_i}(g)(\delta_i+1)}\sum_{J \subset \K^g}
\prod_{j \in J}{{uv-1} \over {(uv)^{\delta_j+1}-1}}E(W \cap D_J^{\circ}/C(g,J))
$$
where $C(g,J)$ is the subgroup of the centralizer of $g$ leaving 
$\bigcap_{j \in J} D_j$ 
invariant.

One shows that for a Kawamata log-terminal $(X_0,D_0)$ pair 
the $E$-function $E(X,D,G)$ of a resolution does not depend 
on the latter but only on $(X_0,D_0,G)$. 
Hence (\ref{efunctionorb}) yields an invariant of Kawamata log-terminal 
$G$-pairs. 
The relation with 
$Ell(X,D,G)$ is the following (cf. \cite{annals}, Prop. 3.14) :

\begin{equation}
{\rm lim}_{\tau \rightarrow i\infty}Ell(X,D,G,z,\tau)=
y^{{-\dim X} \over 2}E(X,D,G,y,1)
\end{equation}
where $y=exp(2 \pi i z)$. 
In particular in non equivariant smooth case the elliptic genus 
for $q \rightarrow 0$ specializes into the Hirzebruch $\chi_y$ genus
(\ref{chisuby}).

On the other hand, in non singular case Hirzebruch (\cite{hirzlevelngenus},
\cite{hirzmodforms}) and Witten \cite{witten87} defined elliptic genera
of complex manifolds which are given by modular forms for the 
subgroup $\Gamma_0(n)$ on level $n$ in $SL_2(\ZZ)$ 
provided the canonical class 
of the manifold in question is divisible by $n$. 
These genera are of course combinations of Chern numbers but for $n=2$ 
one obtains a combination of Pontryagin classes i.e. 
the invariant which depends only on underlying smooth rather 
than (almost) complex structure (this genus was first introduced by S.Ochanine
 cf. \cite{Ochanine} and section \ref{sectionrealsing}). 
These level $n$ elliptic genera 
up to a dimensional factors coincide with specialization 
$z={{\alpha\tau+b} \over n}, \ \ \alpha, \beta \in \ZZ$ for 
appropriate $\alpha, \beta$ specifying particular Hirzebruch 
level $n$ elliptic genus (cf. \cite{invent}, Prop. 3.4).

\subsection{Application I: McKay correspondence for elliptic genus.}
\label{applone}
The classical McKay correspondence is a relation between 
the representations of the  binary dihedral groups 
$G \subset SU(2)$ 
(which are classified according to the root systems 
of type $A_n,D_n, E_6,E_7,E_8$) and the irreducible components of the 
exceptional set of the minimal resolution of $\CC^2/G$. In particular, 
the number of conjugacy classes in $G$ is the same as the number 
of irreducible components of the minimal resolution. The latter is a special 
case of the relation between the Euler characteristic $e(\widetilde {X/G})$ 
of a crepant resolution 
of the quotient $X/G$ of a complex manifold $X$ by an action of a finite 
group $G$ and the data of the action on $X$:

\begin{equation}\label{orbifoldeuler}
 e(\widetilde {X/G})=\sum_{g,h, gh=hg} e(X^{g,h})
\end{equation}

A refinement of the relation (\ref{orbifoldeuler}) for the 
Hodge numbers and motives is given in \cite{bat}, \cite{loeser} 
\cite{looijenga}. In the case when $X$
is projective one has refinement in which the Euler characteristic 
of manifold in
(\ref{orbifoldeuler}) is replaced by elliptic genus of Kawamata log-terminal 
pairs. In fact, more generally, one has the following push forward formula:

\begin{theorem}\label{mckay}
Let $(X;D_X)$ be a Kawamata log-terminal pair which is invariant under 
an effective action of a finite group $G$ on $X$. Let 
$\psi\colon X\to X/G$
be the quotient morphism. Let  $(X/G;D_{X/G})$ be the quotient pair 
in the sense that $D_{X/G}$ is the unique divisor on $X/G$ such that 
$\psi^*(K_{X/G}+D_{X/G})=K_X+D_X$ (cf. \cite{annals}. Def. 2.7)
Then 
$$
\psi_* {\cal Ell}_{orb}(X,D_X,G;z,\tau)={\cal Ell}(X/G,D_{X/G};z,\tau).
$$
\end{theorem}
In particular, for the components of degree zero one obtains:
\begin{equation}
Ell_{orb}(X,D_X,G,z,\tau)=Ell(X/G,D_{X/G},z,\tau)
\end{equation}
In the case when $X$ is non-singular and $X/G$ admits a crepant resolution 
$\widetilde{X/G} \rightarrow X/G$,
for $q=0$ one obtains $\chi_y(\widetilde{X/G})={\chi_y}^{orb}(X,G)$ and hence
for $y=1$ one recovers (\ref{orbifoldeuler}).

\subsection{Application II:Higher elliptic genera and K-equivalences}
\label{appltwo}

Other applications of the push forward formula (\ref{orbblowup}) 
is invariance of higher elliptic genus for the K-equivalences. A question 
posed in (\cite{Rosenberg}) (and answered in \cite{Weinberger}) 
concerns the higher arithmetic genus $\chi_{\alpha}(X)$
of a complex manifold $X$ corresponding to 
a cohomology class $\alpha \in H^*(\pi_1(X),\QQ)$ and defined as 
\begin{equation}\label{higherarithmetic}
  \int_X Td_X \cup f^*(\alpha)
\end{equation} 
where $f: X \rightarrow B(\pi_1(X))$ is 
the classifying map from $X$ to the classifying 
space of the fundamental group of $X$. It asks if  
the higher arithmetic genus $\chi_{\alpha}(X)$ is a birational invariant.
This question is motivated by the Novikov's conjecture: the higher signatures
(i.e. the invariant defined for topological manifold $X$ by 
(\ref{higherarithmetic}) in which the Todd class is replaced
by the $L$-class) are homotopy invariant (cf. \cite{davis}).
Higher $\chi_y$-genus defined by (\ref{higherarithmetic}) with the Todd class
replaced by Hirzebruch's $\chi_y$ class (\cite{hirztop}) comes into 
the correction terms describing the non-multiplicativity of the $\chi_y$ 
in topologically locally trivial 
fibrations $\pi: E \rightarrow B$ of projective manifolds
with non trivial action of $\pi_1(B)$ on the cohomology of the fibers of $\pi$ 
(cf. \cite{maxim} for details).

Recall that two manifolds $X_1,X_2$ are called $K$-equivalent if 
there is a smooth manifold $\tilde X$ and a diagram:

\begin{equation}
\label{kequivalencepicture}
 \begin{matrix} &  & \tilde X & & & \cr
                         & \phi_1\swarrow &  & \searrow \phi_2 & & \cr
                         X_1 & & & & X_2 \cr
\end{matrix}   
\end{equation}
in which $\phi_1$ and $\phi_2$ are birational  
morphisms and $\phi_1^*(K_{X_1})$ and $\phi_2^*(K_{X_2})$ 
are linearly equivalent.

The push forward formula (\ref{orbblowup})
leads to the following:

\begin{theorem}\label{kequivalence} For any $\alpha \in H^*(B\pi,\QQ)$ 
the higher elliptic genus $(\Ell(X) \cup f^*(\alpha),[X])$
is an invariant of $K$-equivalence. Moreover, if $(X,D,G)$ 
and  $(\hat X, \hat D, G)$ 
are $G$-normal and Kawamata log-terminal and if 
$\phi: (\hat X,\hat D) \rightarrow (X,D)$ is $G$-equivariant
such that 
\begin{equation}\label{canonicalclasscondition}
\phi^*(K_X+D)=K_{\hat X}+\hat D
\end{equation}
then 
$$Ell_{\alpha}(\hat X,\hat D,G)=Ell_{\alpha}(X,D,G)$$
In particular the higher elliptic genera (and hence the higher
signatures and $\hat A$-genus)
are invariant for  crepant 
morphisms. The 
specialization into the Todd class
is birationally invariant (i.e. condition of invariance 
(\ref{kequivalencepicture}) is not needed in Todd case).
\end{theorem}

Another consequence is possibility to define higher elliptic genus 
for singular varieties with Kawamata log-terminal singularities 
and for $G$-normal pairs $(X,D)$ (cf. \cite{higherpaper}). 

\subsection{DMVV formula}\label{applthree}
 Elliptic genus comes into a beautiful product 
formula for the generating series for orbifold elliptic genus 
associated with the action of the symmetric group $S_n$ on 
products $X \times ...\times X$
and for which the first case appears in \cite{DMVV} (with 
string-theoretical explanation). A general product 
formula for orbifold elliptic genus of triples is given in \cite{annals}.

\begin{theorem}\label{DMVVpairs}
Let $(X,D)$ be a Kawamata log-terminal pair. For every $n\geq 0$ consider
the quotient of $(X,D)^n$ by the symmetric group $S_n$, which 
we will denote by $(X^n/S_n,D^{(n)}/S_n)$. Here we denote by 
$D^{(n)}$ the sum of pullbacks of $D$ under $n$ canonical projections to $X$. 
Then we have 
\begin{equation}\label{borcherds}
\sum_{n\geq 0}p^n {Ell}(X^n/S_n,D^{(n)}/S_n;z,\tau)=
\prod_{i=1}^\infty \prod_{l,m}
\frac 1 {(1-p^iy^lq^m)^{c(mi,l)}},
\end{equation}
where the elliptic genus of $(X,D)$ is
$$
\sum_{m\geq 0}\sum_l c(m,l)y^lq^m
$$
and $y=\ee^{2\pi\ii z}$, $q=\ee^{2\pi\ii\tau}$.
\end{theorem}

It is amazing that such simple-minded construction as LHS of 
(\ref{borcherds}) leads to the Borcherds lift (cf. \cite{borcherds})
of Jacobi forms.

\subsection{Other applications of elliptic genus} In this section we point out 
other instances in which elliptic genus plays significant role.

\subsubsection{Chiral deRham complex} In work \cite{MSV} for a complex manifold $X$ 
the authors construct a (bi)-graded sheaf $\Omega^{ch}_X$
 of vertex operator algebras 
(with degrees called fermionic charge and conformal weight)
with the differential $d_{DR}^{ch}$ having fermionic degree 1 
and quasi-isomorphic to the deRham complex of $X$. An alternative construction 
in using the formal loop space was given in \cite{KV}. Each component 
of fixed conformal weight has filtration so that graded components 
are:

\begin{equation}
\otimes_{n \ge 1}(\Lambda_{-yq^{n-1}}
T^*_X \otimes \Lambda_{-y^{-1}q^n} T_{X} \otimes S_{q^n}T^*_X \otimes
S_{q^n} T_X)
\end{equation}

In particular it follows that: 

\begin{equation}\label{orbell}
   Ell(X,q,y)=y^{{-\dim X}\over 2}\chi(\Omega^{ch}_X)=
y^{{-\dim X}\over 2}Supertrace_{H^*(\Omega^{ch}_X)} y^{J[0]}q^{L[0]}
\end{equation}
($J[m],L[n]$ are the operators which are part of the vertex algebra structure)).
The chiral complex for orbifolds was constructed in \cite{FrenkelSz}
and the extension of relation (\ref{orbell}) to orbifolds (with discrete torsion)
is discussed in \cite{Matt}.

\subsubsection{Mirror symmetry} Physics definition of mirror 
symmetry in terms of conformal field theory suggests that  
for the elliptic genus, defined as an invariant of a conformal 
field theory, (by an expression similar to 
the last term in (\ref{orbell}), cf. \cite{witten92}) one should have 
for $X$ and its mirror partner $\hat X$ the relation:

\begin{equation}
 Ell(X)=(-1)^{dim X}Ell(\hat X)
\end{equation}
This is indeed the case (cf. \cite{invent}, Remark 6.9)  for the mirror symmetric hypersurfaces 
in toric varieties in the sense of Batyrev .

\subsubsection{Elliptic genus of Landau-Ginzburg models} 
Physics literature (cf. e.g. \cite{Kawai}) also associates 
to a weighted homogeneous polynomial a conformal field theory (Landau-Ginzburg model) 
and in particular the elliptic genus. Moreover it is expected that orbifoldized 
Landau-Ginzburg model will coincide with the conformal filed theory of hypersurface
corresponding to this weighted homogeneous polynomial. 
In particular one expects a certain identity expressing equality of orbifoldized elliptic 
genus corresponding to weighted homogeneous polynomial (or a more general Landau Ginzburg model)
and the elliptic genus of the corresponding hypersurface.
In \cite{MG} the authors 
construct a vertex operator algebra which is related by such type of a correspondence to the 
cohomology of the chiral deRham complex of the hypersurface in $\PP^n$ and in 
particular obtain the expression for the elliptic genus of a hypersurface as
an orbifoldization. Moreover, in \cite{GO} the authors obtain expression for 
the one variable Hirzebruch's genus as an orbifoldization.  

\subsubsection{Concluding remarks} 
There are several other interesting issues which 
should be mentioned in a discussion of elliptic genus. 
It plays important role in work J.Li, K.F.Liu and J.Zhou
(cf. \cite{liu}) in connection with Gopakumar-Vafa conjecture 
(cf. also \cite{GL}).
Elliptic genus was defined for proper schemes with 
1-perfect obstruction theory (\cite{FG}). In fact one has  well 
defined cobordism class in $\Omega^U$ associated to such objects (cf. \cite{kapranov}). 
In the case of surfaces with normal singularities, one can extend 
the above definition of elliptic genus beyond log-terminal singularities
(cf. \cite{waelder3}).
Elliptic genus is central in the study of elliptic cohomology (\cite{Thomas}).
Much of the above discussion of can be extended to equivariant context
(cf. \cite{waelder2}) and a survey of this is given in \cite{waelder1} in this volume.

\section{Quasi-Jacobi forms}

The Eisenstein series 
$e_k(\tau)=\sum_{(m,n) \in \ZZ^2, (m,n) \ne (0,0)} {1 \over {(m \tau+n)^k}}
\ \ (\tau \in \HH)$ 
fails to be modular 
for $k=2$ but the algebra generated by functions $e_k(\tau), k \ge 2$, 
called the algebra of quasi-modular forms on $SL_2(\ZZ)$, 
has many interesting properties
 (cf. \cite{zagier123}). In particular, there is a 
correspondence between quasi-modular forms and real analytic functions 
on $\HH$ which have the same $SL_2(\ZZ)$ transformation 
properties as modular forms. Moreover, the algebra of quasi-modular forms 
has a structure of ${\cal D}$-module and supports an extension 
of Rankin-Cohen operations on modular forms.

In this section we show that there is algebra of functions 
on $\CC \times \HH$ closely related to the algebra of Jacobi forms 
of index zero with similar properties. This algebra is generated by 
Eisenstein series $\sum {1 \over {(z+\omega)^n}}$ (sum over elements
$\omega$ of a lattice $W \subset \CC$).  It has 
description in terms of real analytic functions satisfying 
functional equation of Jacobi forms and having other properties of 
quasi-modular forms mentioned in the last paragraph.
It turns out that the space of functions on $\CC \times \HH$ generated 
by elliptic genera of arbitrary (possibly not
Calabi Yaus) complex manifolds
belong to this algebra of quasi-Jacobi forms.

\bigskip

Recall the following:

\begin{dfn}\label{jacobi}
A weak (resp. meromorphic)
Jacobi form of index $t \in {1 \over 2}\ZZ$ and weight $k$ for a 
finite index subgroup of the Jacobi group
$\Gamma_1^J=SL_2(\ZZ) \propto \ZZ^2$
is a holomorphic 
(resp. meromorphic)
function $\chi$ on $\HH \times \CC$ having expansion
$\sum c_{n,r}q^n\zeta^r$ in $q=exp (2 \pi \sqrt{-1} \tau)$ ($Im\tau$ 
sufficiently large) and
satisfying the following functional 
equations:
$$
\chi({{a\tau+b} \over {c\tau+d}},{z \over {c\tau+d}})=
(c\tau+d)^ke^{{2 \pi i t c z^2} \over {c\tau+d}}\chi(\tau,z)$$ 
$$\chi(\tau,z+\lambda\tau+\mu)=(-1)^{2t(\lambda+\mu)}
e^{-2\pi i t(\lambda^2\tau+2 \lambda z)}\chi(\tau,z)
$$
for all elements 
$[\left( \begin{array}{ccc} a & b \\ c & d \end{array}
\right ), 0 ]$ and $[\left( \begin{array}{ccc} 1 & 0 \\ 0 & 1 \end{array}
\right ), (a,b) ]$ 
in $\Gamma$.
The algebra of Jacobi forms is the bi-graded algebra $J=\oplus J_{t,k}$.
and the algebra of Jacobi forms of index zero is the sub-algebra 
$J_0=\oplus_k J_{0,k} \subset J$. 
\end{dfn} 

For appropriate $l$ a Jacobi form can be expanded in (Fourier)
series in $q^{1 \over l}$ (with $l$ depending on $\Gamma$). We 
shall need below the following real analytic functions: 
\begin{equation}\label{lambdamu}
\lambda(z,\tau)={{z-\bar z} \over {\tau-\bar \tau}}, 
\ \ \ \ \mu(\tau)=
{1 \over {\tau -\bar \tau}}
\end{equation}

They have the following transformation properties:

\begin{equation}
 \lambda({z \over {c\tau+d}},{{a \tau+b} \over {c \tau+d}})=
(c\tau+d)\lambda(z,\tau)-2icz
\end{equation}
$$\lambda(z+m\tau+n,\tau)=\lambda(z,\tau)+m$$
\begin{equation}
\mu({{a \tau+b} \over {c \tau+d}})=(c\tau+d)^2\mu(\tau)
-2ic(c\tau+d)
\end{equation}

\begin{dfn} {\it Almost meromorphic Jacobi form} of weight $k$, index zero
and  depth $(s,t)$  is a (real) meromorphic function in 
$\CC\{q^{1 \over l},z\}[z^{-1},\lambda, \mu]$, with $\lambda, \mu$ given by 
(\ref{lambdamu})  which 

a) satisfies the functional 
equations (\ref{jacobi}) of Jacobi forms of weight $k$ and index zero 
and 

b) which 
has degree at most $s$ in $\lambda$ and at most $t$ in $\mu$. 
\end{dfn}

\begin{dfn}
{\it A quasi-Jacobi form} is a constant 
term of an almost meromorphic 
Jacobi form of index zero 
considered as a polynomial in the functions  $\lambda, \mu$ i.e.
a meromorphic function $f_0$ on $\HH \times \CC$ such that exist meromorphic functions
$f_{i,j}$ such that $f_0+\sum f_{i,j}\lambda^i\mu^j$ is almost meromorphic 
Jacobi form.
\end{dfn}

From algebraic independence of $\lambda, \mu$ over the field 
of meromorphic functions in $q,z$ one deduces:

\begin{prop} $F$ is a quasi-Jacobi of depth $(s,t)$ if and only if:
$$(c\tau+d)^{-k}f({{a\tau+b} \over {c\tau+d}},{z\over {c\tau+d}})=
\sum_{i \le s,j \le t} 
S_{i,j}(f)(\tau,z)({{cz} \over {c\tau+d}})^i({c \over {c\tau+d}})^j$$
$$f(\tau,z+a\tau+b)=\sum_{i \le s} T_i(f)(\tau,z)a^i$$
\end{prop}

To see some basic examples of quasi-Jacobi forms let
us consider the following sequence of functions on $\HH \times \CC$.

\begin{dfn}\label{eisensteindef}
(cf.\cite{Weil}\footnote{These series were recently used
in \cite{Gaberdiel} under the name twisted Eisenstein series.})
$E_n(z,\tau)=\sum_{(a,b) \in \ZZ^2}{1 \over {(z+a\tau+b)^n}}$
\end{dfn}
Such series $E_n(z,\tau)$ 
converges absolutely for $n \ge 3$ and for $n=1,2$ defined
via ``Eisenstein summation'' as $${\sum}_e=lim_{A \rightarrow \infty} 
\sum_{a=-A}^{a=A} (lim_{B \rightarrow \infty} \sum_{b=-B}^{b=B})$$
though we shall omit the subscript $e$.
The series $E_2(z,\tau)$ is related to the Weierstrass function
as follows:
 $$\wp(z,\tau)={1 \over z^2}+
\sum_{(a,b) \in \ZZ, (a,b) \ne 0}{1 \over {(z+a \tau+b)^2}}-
{1 \over {(a \tau+b)^2}}=$$ 
$$=E_2(z,\tau)-lim_{z \rightarrow 0} (E_2(z,\tau)-{1 \over {z^2}})$$
Moreover 
$$e_n=lim_{z \rightarrow 0}(E_n(z,\tau)-{1 \over {z^n}})=
\sum_{(a,b) \in \ZZ, (a,b) \ne 0}{1 \over {(a \tau+b)^n}}$$ is 
the Eisenstein series (notation of \cite{Weil}).
 The algebra of functions of $\HH$ 
generated by the Eisenstein series $e_n(\tau), n \ge 2$ 
 is the algebra of quasi-modular forms 
for $SL_2(\ZZ)$
(cf. \cite{ZagierLandweb}, \cite{zagier123}). 

Now we shall describe the algebra of quasi-Jacobi forms for Jacobi group
$\Gamma^J_1$. We have the following:

\begin{prop}
The functions $E_n$ are weak meromorphic 
Jacobi forms of index zero and weight $n$
for $n \ge 3$.
$E_1$ is a quasi-Jacobi form of index 0 weight 1 and depth $(1,0)$. 
$E_2-e_2$ is a weak Jacobi form of index zero and weight $2$
and $E_2$ is a quasi-Jacobi form of weight 2, index zero and depth $(0,1)$.
\end{prop}

\begin{proof} The first part is due to the absolute convergence 
of series (\ref{eisensteindef}) for $n \ge 3$. We have the following
transformation formulas 
\begin{equation}\label{E1}
E_1({{a\tau+b} \over {c\tau+d}},{z\over {c\tau+d}})=
(c\tau+d)E_1(\tau,z)+{{\pi i c} \over 2}z
\end{equation}
$$E_1(\tau,z+m\tau+n)=E_1(\tau,z)-2 \pi  i m$$
and
\begin{equation}\label{E2}
E_2({{a\tau+b} \over {c\tau+d}},{z\over {c\tau+d}})=
(c\tau+d)^2E_2(\tau,z)-{1 \over 2}\pi i c(c\tau+d)
\end{equation}
$$E_2(\tau,z+a\tau+b)=E_2(\tau,z)$$

First equality in (\ref{E1}) (resp. (\ref{E2})) follows from 
$E_1(\tau,z)={1 \over z}-\sum e_{2k}(\tau)z^{2k-1}$
(resp. 
 $E_2(\tau,z)=
{1 \over z^2}+\sum_k(2k-1)e_{2k}z^{2k-2}$ (cf. \cite{Weil} Ch.3 (10))).
The second in (\ref{E2}) 
is immediate form the definition of Eisenstein summation while the 
second equality in (\ref{E1}) follows from \cite{Weil}.
\end{proof}

\begin{remark} Eisenstein series $e_k(\tau), k \ge 4$ belong to the algebra
of quasi-Jacobi forms. Indeed, one has the following 
identities (cf. \cite{Weil}, (7),(35) in ch IV):
$$E_4=(E_2-e_2)^2-5e_4; \ \ \ E_3^2=(E_2-e_2)^2-15e_4(E_2-e_2)-35e_4$$
\end{remark}

\begin{prop} The algebra of Jacobi forms (for $\Gamma^J_1$) 
of index zero and weight $t$
($t \ge 2$) is generated by $E_2-e_2,E_3,E_4$.
\end{prop}

A short way to show this is to notice that 
the ring of such Jacobi forms is isomorphic 
to the ring of cobordisms of $SU$-manifolds
 modulo flops (cf. section \ref{start})
 via isomorphism sending a complex manifold $X$ of 
dimension $d$ to $Ell(X) \cdot ({{\theta'(0)} \over {\theta(z)}})^d$.
This ring of cobordisms  
in turn is isomorphic to $\CC[x_1,x_2,x_3]$ where $x_1$ is cobordism 
class of a K3 surface and $x_2,x_3$ are cobordism classes of certain 
four-  and 
six-manifolds
(cf. \cite{Totaro}).
The graded algebra $\CC[E_2-e_2,E_3,E_4]$ is isomorphic to the 
same ring of polynomials (cf. examples \ref{examples})
and the claim follows.

\begin{prop} The algebra of quasi-Jacobi forms is the algebra of functions
on $\HH \times \CC$ generated by functions $E_n(z,\tau)$
and $e_2(\tau)$.
\end{prop}

\begin{proof} First notice that the coefficient of $\lambda^s$ for 
an almost meromorphic 
Jacobi form $F(\tau,z)=\sum_{i \le s} f_i\lambda^i$ 
of depth $(s,0)$ is holomorphic Jacobi form of index zero and weight
$k-s$ i.e. by above proposition is a polynomial in $E_2-e_2,E_3,..$.
Moreover $f_0-E_1^sf_s$ is a quasi-Jacobi form of index zero and 
weight at most $s-1$. Hence by induction the ring of quasi-Jacobi forms
of index zero and depth $(*,0)$ can be identified with 
$\CC[E_1,E_2-e_2,E_3,...]$. Similarly, the coefficient $\mu^t$ 
of an almost meromorphic 
Jacobi form $F=\sum_{j \le t}(\sum f_{i,j}\lambda^i)\mu^j$ is an almost 
meromorphic Jacobi form of depth $(s,0)$ and moreover 
$F-(\sum_{i}f){i,s}\lambda^i)E_2^t$ has depth $(s',t')$ with 
$t'<t$. The claim follows.
\end{proof}

Here is an alternative description of the algebra of quasi-Jacobi forms:

\begin{prop}\label{zagiercoefficients}
 The algebra of functions generated by the coefficients 
of the Taylor expansion in $x$ of the function: 
$${{\theta(x+z)\theta'(0)} \over {\theta(x)\theta(z)}}-
({1 \over x}+{1 \over z})=\sum_{i \ge 1} 
F_ix^i$$ is the algebra of quasi-Jacobi forms 
(for $SL_2(\ZZ)$).
\end{prop}

\begin{proof} The transformation formulas for theta function
(cf. \cite{TanMolk}):
\begin{equation}
\theta({{a\tau+b} \over {c\tau+d}},{z\over {c\tau+d}})=
\zeta(c\tau+d)^{1 \over 2}e^{{\pi ic z^2} \over {c\tau+d}}\theta(\tau,z)
\end{equation}
$$\theta'({{a\tau+b} \over {c\tau+d}},0)=\zeta(c\tau+d)^{3\over 2}
\theta'(\tau,0)$$
$$\theta(\tau, z+m\tau+n)=
(-1)^{m+n}e^{-2\pi i mz-\pi i m^2 \tau}\theta(\tau,z)$$
imply that  
\begin{equation}\Phi(x,z,\tau)=
 {{x\theta(x+z)\theta'(0)} \over {\theta(x)\theta(z)}}
\end{equation} 
satisfies the 
following functional equations:
\begin{equation}\label{transformationsPhi}
\Phi({{a\tau+b} \over {c\tau+d}},{x\over {c\tau+d}},
{z\over {c\tau+d}})=
%{1 \over {c\tau+d}}
e^{{2 \pi i c zx} \over {c\tau+d}}\Phi(x,z,\tau)
\end{equation}
$$\Phi(x,z+m\tau+n,\tau)=e^{2 \pi i m x}\Phi(x,z,\tau)$$

In particular in expansion 
\begin{equation}{{d^2(log(\Phi))} \over {dx^2}}=\sum H_ix^i
\end{equation} 
the left hand side is invariant under transformations in 
(\ref{transformationsPhi}) and 
the coefficient $H_i$ is a Jacobi form of weight $i$ and index zero for any 
$i$.
Moreover the coefficients $F_i$ in $\Phi(x,z,\tau)=1+\sum F_i(z,\tau)x^i$
are polynomials in $F_1$ and $H_i$. What remains to show is that 
$E_i$'s determine $F_1,H_i, i \ge 1$ and vice versa.
Recall that $E_i$ has index zero (invariant with respect to shifts) and 
weight $i$. 
We shall use the following expressions:
\begin{equation}\label{zagier}
\Phi(x,z,\tau)={{x+z}\over {z}}exp(\sum_{k>0}
{2 \over {k!}}[x^k+z^k-(x+z)^k]G_k(\tau))
\end{equation}
where 
\begin{equation}
 G_k(\tau)=-{{B_k} \over {2k}}+\sum_{l=1}^{\infty}\sum_{d \vert l}(d^{k-1})q^l
\end{equation}
(cf. \cite{Zagier}). On the other hand one has:
\begin{equation}
E_n(z,\tau)={1 \over {z^n}}+(-1)^n\sum_{2m \ge n}^{\infty}{{2m-1} \choose {n-1}}
e_{2m}z^{2m-n}
\end{equation}
(cf. \cite{Weil} Ch.III.sec.7 (10))
where 
\begin{equation}\label{e2m}
e_{2m}={\sum}'({1 \over {m \tau+n}})^{2m}={{2(2\pi \sqrt{-1})^{k}}\over{(k-1)!}}G_k \ \ 
({\rm for } \ k=2m)
\end{equation}
(cf. \cite{Weil} Ch.III, sect 7 and \cite{ZagierLandweb} p.220).
We have:
\begin{equation}
{{d^2log(\Phi(x,z,\tau))} \over {dx^2}}=
\sum_{i \ge 1}{{(-1)^iix^{i-1}} \over 
{z^{i+1}}}+
\sum_{i \ge 2}{2 \over {(i-2)!}}[x^{i-2}-(x+z)^{i-2}]G_i(\tau)
\end{equation}
The coefficient of $x^{l-2}$ for $l \ge 2$ in Laurent expansion hence 
is:
\begin{equation}
{{(-1)^{l-1}(l-1)}\over{z^{l}}}-
\sum_{i \ge 2, i > l}{2 \over {(i-2)!}}{{i-2} \choose {l-2}}z^{i-l}G_i(\tau)=
\end{equation}
$${{(-1)^{l-1}(l-1)}\over{z^{l-1}}}-
-\sum_{i \ge 2, i >l}
{1 \over {(2 \pi \sqrt{-1})^i}}(i-1){{i-2} \choose {l-2}}z^{i-l}
e_i=
$$
$${{(-1)^{l-1}(l-1)}\over{z^l}}-(l-1){1 \over {(2 \pi \sqrt{-1})^l}}
\sum_{i \ge 2, i>l}{{i-1} \choose {l-1}}
e_i({z \over {2 \pi \sqrt{-1}}})^{i-l} 
$$
(using (\ref{e2m}) and identities with binomial coefficients). This yields
$$H_{l-2}({2 \pi \sqrt{-1}}z,\tau)
=(-1)^{l-1}{{(l-1)} \over {({2 \pi \sqrt{-1}})^l}}
(E_l-e_l)$$ and the claim follows since 
(15) in \cite{Zagier} yields:

\begin{equation}
 F_1(z,\tau)={1 \over z}-2\sum_{r \ge 0} G_{r+1}{{z^r} \over {r!}}=
{1 \over z}- {1 \over {(2 \pi \sqrt{-1})}} \sum_{r \ge 0}
e_r({z \over 
{2 \pi \sqrt{-1}}})^r
\end{equation}
i.e.
\begin{equation}
F_1(2 \pi i \sqrt{-1} z,\tau)
={1 \over {2 \pi i \sqrt{-1}}}E_1(z, \tau)
\end{equation} 

\end{proof}

\begin{remark} The algebra of quasi-Jacobi forms 
$\CC[e_2,E_1,E_2,...]$
in  closed under
differentiations $\partial_{\tau}, \partial_z$. 
Indeed one has: 
$$2 \pi i {\partial E_1 \over {\partial \tau}}=E_3-E_1E_2. \ \ \ \ \
     {\partial E_1 \over {\partial z}}=-E_2$$

$$2 \pi i {\partial E_2 \over {\partial \tau}}=3E_4-2E_1E_3-E_2^2.
\ \ \ \ \      {\partial E_2 \over {\partial z}}=-2E_3$$
and hence $\CC[..E_i..]$ is a $\cal D$-module where
$\cal D$ be the ring of differential operators generated by 
over the ring of holomorphic Jacobi group invariant functions on 
$H \times \CC$ by $\partial \over {\partial \tau}$
and $\partial \over {\partial z}$.
As is clear from the above discussion, 
the ring of Eisenstein series $\CC[...E_i...]$ has natural 
identification with the ring of real valued almost meromorphic 
Jacobi forms 
$\CC[E_1^*,E_2^*,E_3,....]$ 
on $\HH \times \CC$ having index zero where
\begin{equation}
E_1^*=E_1+2 \pi i{{Im x} \over {Im \tau}}, E_2^*=E_2+{1 \over {Im \tau}} 
\end{equation}
\end{remark}

\begin{theorem} The algebra of quasi-Jacobi forms of depth $(k,0), k \ge 0$
is isomorphic to the
algebra of complex unitary cobordisms modulo flops.
\end{theorem}

In another direction the depth of quasi-Jacobi forms allowing to ``measure'' 
the deviation of elliptic genus of a non-Calabi Yau 
manifold from being Jacobi form.

\begin{theorem} Elliptic genera of manifolds of dimension at most $d$ span 
the subspace of forms of depth $(d,0)$ in the algebra of quasi-Jacobi forms.
If a complex manifold satisfies $c_1^k=0, c_1^{k-1} \ne 0$ 
\footnote{more generally, $k$ is the smallest among 
$i$ with $c_1^i \in Ann(c_2,..c_{dim M})$; an example of such manifold is 
an $n$-manifold having $n-k$-dimensional Calabi Yau factor.}
then 
its elliptic genus is a quasi-Jacobi form of depth $(s,0)$ where $s \le k-1$.
\end{theorem}

\begin{proof}
 It follows from the proof of proposition (\ref{zagiercoefficients})
that $${{d^2log \Phi} \over {dx^2}}=\sum_{i \ge 2} (-1)^{i-1}
{{i-1} \over {(2 \pi \sqrt{-1})^i }}(E_i-e_i)x^{i-2}$$
which yields: 
\begin{equation}
\Phi=e^{E_1x}\Prod_i e^{{1 \over i}(-1)^{i-1}
{{i-1} \over {(2 \pi \sqrt{-1})^{i} }}(E_i-e_i)x^{i}}
\end{equation}
The Hirzebruch's characteristic series is (cf. (\ref{ellipticclass})):
$$\Phi({x \over {2 \pi i}})({{\theta(z)} \over {\theta'(0)}})$$
Hence if $c(TX)=\Pi(1+x_k)$ then 
\begin{equation}\label{charactseries}
Ell(X)=
({\theta(z) \over {\theta'(0)}})^{dim X}
\Prod_{i,k} e^{E_1x_k}e^{{(-1)^{i-1}{(i-1)} \over i}
(E_i-e_i)x^i_k}[X]=
\end{equation}
$$({\theta(z) \over {\theta'(0)}})^{dim X}
e^{c_1(X)E_1}\Prod_{i,k}e^{{(-1)^{i-1}{(i-1)} \over i}(E_i-e_i)x^i_k}[X]
$$
(here $[X]$ is the fundamental class of $X$) i.e. if $c_1=0$ 
elliptic class is polynomial in $E_i-e_i$ with $i \ge 2$ and hence
elliptic genus is a Jacobi form (cf. \cite{krichever}).
Moreover of $c_1^k=0$ then the degree of this polynomial 
is at most $k$ in $E_1$ and the claim follows.
\end{proof}

\begin{example}\label{examples} Expression (\ref{charactseries}) can be used to get 
formulas for the elliptic genus of specific examples in terms of 
Eisenstein series $E_n$. For example for a surface in $\PP^3$ 
having degree $d$ one has
$$ (E_1^2({1 \over 2}d^2-4d+8)d+(E_2-e_2)({d^2 \over 2}-2)d)({\theta(z) \over
\theta'(0)}^2)
$$ 
In particular for $d=1$ one obtains:
$$({9 \over 2}E_1^2-{3 \over 2}(E_2-e_2))({\theta(z) \over
\theta'(0)})^2
$$
One can compare this with the double series which is a special case
of the general formula for 
elliptic genus of toric varieties in \cite{invent}. This lead to a two 
variable version of the identity discussed in Remark 5.9 in \cite{invent}.
In fact following \cite{gunnels} one can define sub-algebra 
``toric quasi-Jacobi forms'' of the algebra of quasi-Jacobi forms
extending toric quasi-modular forms considered in \cite{gunnels}. This 
issue will be addressed elsewhere.
\end{example}

Next let us consider one more similarity between meromorphic Jacobi 
forms and modular forms: there is a natural non commutative 
deformation of the ordinary
product of Jacobi forms similar to the deformation of the product 
modular forms constructed using Rankin-Cohen brackets (cf. \cite{zagier123}).
In fact we have the following Jacobi counterpart of the Rankin-Cohen
brackets:

\begin{prop}\label{bracket}
 Let $f$ (resp. $g$) be a Jacobi form of index zero and weight $k$ 
(resp. $l$). Then $$[f,g]={k}(\partial_{\tau}f-{1 \over {2 \pi i }}E_1\partial_zf)g-
{l}(\partial_{\tau}g-{1 \over {2 \pi i }}E_1(z,\tau)\partial_zg)f$$
is a Jacobi form of weight $k+l+2$.
More generally, let $$D=\partial_{\tau}-{1 \over {2 \pi i }}E_1\partial_z$$
Then the Cohen-Kuznetsov series (cf. \cite{zagier123}):
$$\tilde f_D(z, \tau, X)
=\sum_{n=0}^{\infty} {{D^nf(z,\tau) X^n} \over 
{n!(k)_n}} $$
(here $(k)_n=k(k+1)...(k+n-1)$ is the Pochhammer symbol) satisfies:
$$\tilde f_D({{a \tau+b} \over {c\tau+d}},{z \over {c\tau+d}},{z \over 
{c\tau+d}},{X \over {(c\tau+d)^2}})=(c\tau+d)^kexp({c \over {c \tau+d}}
{X \over 2 \pi i})f_D(\tau,z,X)
$$
$$\tilde f_D(\tau,z+a\tau+b,X)=\tilde f_D(\tau, z,X)$$
In particular the coefficient 
$${{[f,g]_n} \over {(k)_n(l)_n}}$$
of $X^n$ in $\tilde f_D(\tau,z,-X)\tilde g_D(\tau,z,X)$ is a Jacobi form of
weight $k+l+2n$. It is given explicitly in terms of $D^{i}f,D^{j}g$ by 
the same formulas as the classical RC brackets.
\end{prop}

\begin{proof} The main point is that the operator 
$\partial_{\tau}-{1 \over {2\pi i }}E_1\partial_z$ has the same 
deviation from transforming Jacobi form to Jacobi  
as $\partial_{\tau}$ on modular forms. Indeed:
$$(\partial_{\tau}-{1 \over {2\pi i }}E_1\partial_z)
f({{a \tau+b}\over {c \tau+d}},{z \over {c\tau+d}})=$$
$$\{kc(c\tau+d)^{k+1}f(\tau,z)+zc(c\tau+d)^{k+1}\partial_zf(\tau,z)+
(c\tau+d)^{k+2}\partial_{\tau}f(\tau,z)$$
$$-{1 \over {2 \pi i }}[(c\tau+d)E_1(\tau,z)+2 \pi icz]
(c\tau+d)^{k+1}\partial_zf)\}=
$$
$$(c\tau+d)^{k+2}[f(\tau,z)-{1 \over {2 \pi i }}E_1f_z]+
kc(c\tau+d)^{k+1}f(\tau,z)
$$
Moreover:
$$(\partial_{\tau}-{1 \over {2\pi i }}E_1\partial_z)f(\tau,
z+a\tau+b)=f_{\tau}+af_z-{1 \over {2 \pi i}}(E_1-2\pi i a)f_z=
(\partial_{\tau}-{1 \over {2\pi i }}E_1\partial_z)f(\tau,z)
$$ 
The rest of the proof is the same as in \cite{zagier123}.
\end{proof}

\begin{remark} The brackets introduced in Prop. \ref{bracket}
are different from the Rankin-Cohen bracket introduced in \cite{choie}.
\end{remark}

\section{Real singular varieties}\label{sectionrealsing}

Ochanine genus of an oriented differentiable manifold $X$ can be defined 
using the following 
series with coefficients in ${\bf Q}[[q]]$ 
as the Hirzebruch characteristic power series
(cf. volume \cite{Landweber.Stong} and references there):
\begin{equation}\label{ochanineseries}
Q(x)={{x/2}\over {\sinh( x/2)}} \prod_{n=1}^{\infty}
  [{{(1-q^n)^2} \over {(1-q^n\ee^x)(1-q^n\ee^{-x})}}]^{(-1)^n} 
\end{equation}

As was mentioned in section \ref{relationtoother}, this genus is
a specialization of the two variable elliptic genus (at $z={1\over 2}$). 
Evaluation of Ochanine genus of a manifold 
using (\ref{ochanineseries}) and  
viewing the result as function of $\tau$ on the upper half plane 
(where $q=\ee^{2 \pi \ii \tau}$) yields a 
modular form on $\Gamma_0(2) \subset SL_2(\ZZ)$ 
(cf. \cite{Landweber.Stong}).

In this section we discuss elliptic genera for real algebraic varieties.
It particular we address B.Totaro proposal that  ``it should 
be possible to define Ochanine genus for a large class of compact 
oriented real analytic spaces'' (cf. \cite{Totaro}). 
In this direction \cite{Totaro} contains the following result:

\begin{theorem} Quotient of $MSO$ by ideal generated by 
oriented real flops and complex flops \footnote{i.e. the ideal 
generated by $X'-X$ where $X'$ and $X$ are related by real 
or complex flop} is 
$$\ZZ[\delta, 2\gamma, 2\gamma^2,2 \gamma^4..]$$
with $\CC\PP^2$ (resp. $\CC\PP^4$) corresponding to $\delta$
(resp. $2 \gamma+\delta^2$). This quotient ring is the 
the image of $MSO_*$ under the Ochanine genus.
\end{theorem}

In particular Ochanine genus of a small resolution is independent
of its choice for singular spaces
having singularities only along non-singular strata and having in normal 
directions only singularities which are cones in $\RR^4$ or $\CC^4$.

Our goal is to find a wider class of singular real algebraic varieties 
for which Ochanine genus of a resolution is independent of 
a choice of the latter.

\subsection{Real Singularities} 

For the remainder of this paper ``real algebraic variety'' means an 
{\it oriented} quasi-projective variety $X_{\RR}$
over $\RR$, $X(\RR)$ is the set of its $\RR$-points with Eucledian topology,
$X_{\CC}=X_{\RR} \times_{Spec \RR} Spec \CC$ is the complexification 
and $X(\CC)$ the analytic space of its complex points. We also assume that 
$dim_{\RR} X(\RR)=dim_{\CC} X(\CC)$.

\begin{dfn} A real algebraic variety $X_{\RR}$ as above  
is called $\QQ$-Gorenstein log-terminal
if the analytic space $X(\CC)$ is $\QQ$-Gorenstein log-terminal.
\end{dfn}

\begin{example}\label{firstexample} Affine 
variety 
\begin{equation}\label{equationfirstexample}
x_1^2-x_2^2+x_3^2-x_4^2=0
\end{equation}
in $\RR^4$
is 3-dimensional Gorenstein log-terminal and admit a crepant resolution.
\end{example}

Indeed it is well known that complexification of 
Gorenstein singularity (\ref{equationfirstexample})
admits a small (and hence crepant) resolution 
having $\PP^1$ as its exceptional set.

\begin{example}\label{secondexample} The 3-dimensional 
complex cone in $\CC^4$ given by 
$z_1^2+z_2^2+z_3^2+z_4^2=0$ considered as codimension 
two sub-variety of $\RR^8$ is a $\QQ$-Gorenstein log-terminal variety 
over $\RR$ and its complexification admits a crepant resolution.

Indeed, this codimension two sub-variety is a real analytic 
space which is the intersection of two 
quadrics in $\RR^8$ given by 
\begin{equation}\label{ci}
a_1^2+a_2^2+a_3^2+a_4^2-b_1^2-b_2^2-b_3^2-b_4^2=0
\end{equation}
$$a_1b_1+a_2b_2+a_3b_3+a_4b_4=0$$
(here $a_i=Re z_i, b_i=Im z_i$).
The complexification, is the cone over complete intersection of 
two quadrics in $\PP^7$. Moreover, 
the defining equations of this complete intersection 
after the change of coordinates:
$x_i=a_i+\sqrt{-1}b_i, y_i=a_i-\sqrt{-1}b_i$
become:
\begin{equation}\label{cip}
x_1^2+x_2^2+x_3^2+x_4^2=0
\end{equation}
 $$y_1^2+y_2^2+y_3^2+y_4^2=0$$ 
The singular locus is the union of two disjoint 2-dimensional 
quadrics and singularity along each is $A_1$ (i.e. the 
intersection of the transversal to it in $\PP^7$ has $A_1$-singularity).
To resolve (\ref{ci}), one can blow up $\CC^8$ at the origin which results 
in $\CC$-fibration over complete intersection (\ref{cip}). It 
can be resolved by small resolutions along two non-singular 
components of singular locus of (\ref{cip}).
A direct calculation shows \footnote{e.g. by considering the 
order of the pole of the form 
${{dx_2 \wedge dx_3 \wedge dx_4 \wedge dy_2 \wedge dy_3 
\wedge dy_4} \over {x_1y_1}}$ along the intersection of 
the exceptional locus of blow up $\tilde {\CC^8}$
of $\CC^8$ with the proper preimage of (\ref{cip}) in
 $\tilde {\CC^8}$}
that we have a log-terminal resolution
of Gorenstein singularity which is the complexification of (\ref{ci}).

\end{example}

\subsection{Elliptic genus of resolutions of real varieties
with $\QQ$-Gorenstein log-terminal singularities.}\label{lastsection}

Let $X$ be a real algebraic manifold and let $D=\sum \alpha_kD_k$ $(\alpha_k 
\in \QQ)$ be a divisor on the 
complexification $X_{\CC}$ of $X$ (i.e. $D_k$ are irreducible components 
of $D$). Let $x_i$ denote the Chern roots
of the tangent bundle of $X_{\CC}$ and $d_k$ are the classes corresponding 
to $D_k$ (cf. section \ref{ellipticgenus}). 

\begin{dfn} Let $X$ be a real algebraic manifold and $D$ a divisor on 
complexification $X_{\CC}$ of $X$. 
The Ochanine class $\Ell_{\cal O}(X,D)$ of pair $X,D$ is specialization:
$$\Ell (X_{\CC},D,q,z={1 \over 2})$$
of the two variable elliptic class of pair $\Ell (X_{\CC},D,q,z)$ given by:
\begin{equation}
\Bigl(\prod_l 
\frac{(\frac{x_l}{2\pi\ii})\theta(\frac{x_l}{2\pi\ii}-z)\theta'(0)} 
{\theta(-z)\theta(\frac{x_l}{2\pi\ii})} 
\Bigr)\times 
\Bigl(\prod_k 
\frac{\theta(\frac{d_k}{2\pi\ii}-(\alpha_k+1)z)\theta(-z)} 
{\theta(\frac{d_k}{2\pi\ii}-z)\theta(-(\alpha_k+1)z)} 
\Bigr)
\end{equation}
Ochanine elliptic genus of pair $(X,D)$ as above is
\begin{equation}
Ell(X_{\RR},D)=\sqrt{\E ll(X_{\CC},D,q,{1 \over 2})}
 \cup cl(X({\RR}))[X({\CC})]
\end{equation}   
\end{dfn}

Here $\sqrt{\E ll}$ denotes the class corresponding to 
the unique series with constant term equal to 1 and having $\E ll$ as its 
square. 

The above class of pair is the class (\ref{formulaorbifoldelliptic}) 
considered in definition \ref{orbifoldelliptic} in
section \ref{ellipticgenus} 
 with group $G$ being trivial. One can defined
orbifold version of this class as well specializing 
(\ref{formulaorbifoldelliptic})
to $z={1 \over 2}$. See \cite{Duke} for further discussion of 
the class $\E ll(X,D)$.

The relation with Ochanine's definition is as follows: if $D$ is trivial 
divisor on $X_{\CC}$ then the result coincides with the genus \cite{Ochanine}.
More precisely, we have the following:

\begin{lemma}\label{nonsingularreal}
 Let $X_{\RR}$ be a real algebraic manifold with non singular 
complexification $X_{\CC}$. Then
$$Ell(X_{\RR})=
\sqrt{\E ll(T_{X({\CC})})} \cup  cl(X({\RR}))[X({\CC})]$$
\end{lemma}

\begin{proof} Indeed, we have 
\begin{equation}
0 \rightarrow T_{X({\RR})} \rightarrow T_{X({\CC})}\vert_{X({\RR})}
\rightarrow T_{X({\RR})} \rightarrow 0
\end{equation}
with the identification of the normal bundle to $X_{\RR}$ with 
its tangent bundle given by the multiplication by $\sqrt{-1}$. Hence 
$\E ll(X_{\RR})^2=i^*(\E ll X_{\CC})$ where $i: X_{\RR} \rightarrow
X_{\CC}$ is the canonical embedding. Now the lemma follows from
the identification of the characteristic series (\ref{ochanineseries})
and specialization $z={1 \over 2}$ of the series in (\ref{ellipticclass})
(cf. \cite{invent}) and  
the identify which is just a definition of the class 
$cl_Z \in H^{{\rm dim}_{\RR}Y-
{\rm dim}_{\RR}}Z$ of 
a sub-manifold $Z$ of a manifold $Y$: $cl_Z \cup \alpha [Y]=
i^*(\alpha) \cap [Z]$ for any $\alpha \in H^{{\rm dim}_{\RR}Z}(Y)$.
Indeed, we have the following:
\begin{equation}
Ell(X_{\RR})=\E ll (T_{X({\RR})})[X({\RR})]=
\sqrt{\E ll(T_{X({\CC})})} \vert_{X({\RR})})[X({\RR})]=
\end{equation}
$$\sqrt{\E ll(T_{X({\CC})})} \cup  cl(X({\RR}))[X({\CC})]
$$

\end{proof}

Our main result in this section is the following:

\begin{theorem}\label{realsing} Let $\pi: (\tilde X,\tilde D) 
\rightarrow (X,D)$ be a resolution 
of singularities 
of a real algebraic pair with $\QQ$-Gorenstein log-terminal 
singularities i.e. $K_{\tilde X}+\tilde D=\pi^*(K_X+D)$. 
Then the elliptic genus of the pair $(\tilde X,\tilde D)$
is independent of resolution.
In particular if real algebraic variety $X$ has a crepant resolution then 
its elliptic genus is independent of a choice of crepant resolution.
\end{theorem}

\begin{proof}
 Indeed for a blow up $f: (\tilde X,\tilde D) \rightarrow 
(X,D)$ we have:   
\begin{equation} f_*(\sqrt{\Ell (\tilde X,\tilde D,q,{1 \over 2}})=
\sqrt{\Ell(X,D,q,{1 \over2})} 
\end{equation}
This is a special case of the push-forward formula (\ref{orbblowupformula})
in theorem \ref{orbblowup} (with $G$ being trivial group).
Hence
\begin{equation} 
\E ll_{\cal O}(X_{\RR},D)=\sqrt{\E ll(X_{\CC},D,q,{1 \over 2})}
 \cup cl(X_{\RR})[X_{\CC}]=
\end{equation}
$$\sqrt{\Ell (\tilde X_{\CC},\tilde D,q,{1 \over 2}}) \cup 
f^*([X_{\RR}] \cap [X_{\CC}])=\Ell (\tilde X_{\RR}, \tilde D)  
$$
as follows from projection formula
since $f^*(cl[X_{\RR}])=[cl \tilde X_{\RR}]$
and since $f_*$ is identity on $H_0$.

For crepant resolution one has $D=0$ and hence by lemma \ref{nonsingularreal}
the elliptic genus of $X_{\RR}$ is the Ochanine genus of the real 
manifold which is its crepant resolution.  

\end{proof}

\begin{remark}Examples \ref{firstexample} and \ref{secondexample}
show that singularities admitting crepant resolution include
real 3-dimensional cones and real points of complex 3-dimensional 
cones.
\end{remark}

%\section{Loop spaces and Vertex Operator algebras}

%\section{Gopakumar-Vafa and Nekrasov conjectures}

\end{document}